\documentclass[12pt,oneside,reqno]{amsart}
\usepackage{amssymb}
\usepackage{}
\usepackage{txfonts}
\usepackage{bbm}
\usepackage{amsmath}
\usepackage{graphicx}
\usepackage{mathrsfs}
\usepackage{stmaryrd}
\usepackage{amsfonts}
\usepackage{enumerate,amsmath,amssymb,amsthm}
\usepackage[hypertexnames=false,bookmarksopen=true,linktocpage=true,pdfstartview={XYZ null null 1.25}]{hyperref}
\usepackage{color}

\numberwithin{equation}{section}

\newcommand{\be}{\begin{eqnarray}}
\newcommand{\ee}{\end{eqnarray}}
\newcommand{\ce}{\begin{eqnarray*}}
\newcommand{\de}{\end{eqnarray*}}
\newtheorem{theorem}{Theorem}[section]
\newtheorem{lemma}[theorem]{Lemma}
\newtheorem{remark}[theorem]{Remark}
\newtheorem{definition}[theorem]{Definition}
\newtheorem{proposition}[theorem]{Proposition}
\newtheorem{Examples}[theorem]{Example}
\newtheorem{corollary}[theorem]{Corollary}

\def\[{{\Big[}}
\def\]{{\Big]}}
\def\<{{\langle}}
\def\>{{\rangle}}
\def\({{\Big(}}
\def\){{\Big)}}

\def\bx{{\mathbf{x}}}

\def\dif{{\mathord{{\rm d}}}}

\def\min{{\mathord{{\rm min}}}}

\def\={&\!\!=\!\!&}
\def\bt{\begin{theorem}}
\def\et{\end{theorem}}
\def\bl{\begin{lemma}}
\def\el{\end{lemma}}
\def\br{\begin{remark}}
\def\er{\end{remark}}

\def\bd{\begin{definition}}
\def\ed{\end{definition}}
\def\bp{\begin{proposition}}
\def\ep{\end{proposition}}
\def\bc{\begin{corollary}}
\def\ec{\end{corollary}}
\def\bx{\begin{Examples}}
\def\ex{\end{Examples}}

\def\cE{{\mathcal E}}

\def\cH{{\mathcal H}}

\def\mE{{\mathbb E}}

\def\mR{{\mathbb R}}

\def\geq{\geqslant}
\def\leq{\leqslant}

\allowdisplaybreaks

\title{\bf{The G-Martingale Approach for G-Utility Maximization}}
\author{
\ Qiguan Chen$^{1}$ Yulin Song$^{1}$\   Zengwu Wang$^{2}$ \ Zengting Yuan$^{2}$\\
{\em $^{1}$D\MakeLowercase {epartment} \MakeLowercase{of} M\MakeLowercase{athematics}, N\MakeLowercase{anjing}
U\MakeLowercase{niversity}, N\MakeLowercase{anjing}, {\rm 210093}, C\MakeLowercase {hina}}\\
{\em $^{2}$I\MakeLowercase {nstitute} \MakeLowercase{of} F\MakeLowercase{inance} \MakeLowercase{and} B\MakeLowercase{anking}, C\MakeLowercase{hinese} A\MakeLowercase{cademy} \MakeLowercase{of}
S\MakeLowercase{ocial} S\MakeLowercase{ciences}, B\MakeLowercase{eijing}, {\rm 100710}, C\MakeLowercase {hina}}
}
\begin{document}
\footnotetext[1]{Qiguan Chen: 171850680@smail.nju.edu.cn }
\footnotetext[2]{Yulin Song: ylsong@nju.edu.cn, supported by National Natural Science Foundation of China (No.11971227, 11790272)}
\footnotetext[3]{Zengwu Wang: the corresponding author, zwwang@amss.ac.cn}
\footnotetext[4]{Zengting Yuan: yztifb@cass.org.cn}
\maketitle

\begin{abstract}
\noindent In this paper, we study representative investor's G-utility maximization problem by G-martingale approach in the framework of G-expectation space proposed by Peng \cite{Pe19}. Financial market has only a bond and a stock with uncertainty characterized by G-Brownian motions. The routine idea of \cite{Wxz} fails because that the quadratic variation process of a G-Brownian motion is also a stochastic process. To overcome this difficulty, an extended nonlinear expectation should be pulled in. A sufficient condition of G-utility maximization is presented firstly. In the case of log-utility, an explicit solution of optimal strategy can be obtained by constructing and solving a couple of G-FBSDEs, then verifying the optimal strategy to meet the sufficient condition. As an application, an explicit solution of a stochastic interest model is obtained by the same approach. All economic meanings of optimal strategies are consistent with our intuitions.

\vskip0.5cm \noindent{\bf Keywords:} G-Expectation, G-Brownian Motion, G-Martingale, G-Utility.\vspace{1mm}\\
\noindent{{\bf JEL Classification} D81  G11  G12}
\end{abstract}

\section{Introduction}

Research on the distinction between risk and uncertainty has lasted more than one hundred years. In 1921, almost at the same time, Keynes \cite{Keynes} and Knight \cite{Kni} affirmed that risk and uncertainty are different conceptions in economic environment. Simply put, risk means that  investors (or decision makers) know all states of the economic environment and their corresponding probabilities, but can not conclude which state will occur later. This can be characterized well by probability theory. However, uncertainty means that  investors know states of economic environment, but can not identify probabilities of these states. For instance, tossing a fair coin, decision makers know that there are only two states, head and tail, and the probability of each state is a half. By contrast, for tossing an unfair coin with the unknown degree of its unfairness, decision makers still know both states, but can not claim their probabilities. In this case, probability of each state is not a constant, but an interval probably containing a half.

It is known that probability theory is just suitable for dealing with risk. However, how to deal with uncertainty is still an open problem up to now. Keynes \cite{Keynes} firstly attempted to redefine probability theory in the view of logics. Afterwards, there are two main routes in this research field. One is extending probability theory, that is, replacing additivity in the axiomatic characterization of probability by monotonicity of set functions. Main results of this route are Choquet capacity and Choquet integration, put forward by Choquet \cite{Cho} and his following researchers. Distorted probability is another remarkable notion in between probability and Choquet capacity. It is a tractable composite function of an increasing function and a probability measure. So it has better applications in finance and insurance, cf. Wang \cite{Wan} and the followed. The other route is nonlinear methods based on the max-min probability. Gilboa and Schmeidler \cite{GS} presented the axiomatic characterization of max-min expected utilities, which defines a new kind of Choquet capacity indeed. 

However, Bayesian rule of both Choquet capacity and max-min probability can not obtain unanimous approvals. This flaw restricts their applications in economics, finance and other fields, especially for their related dynamic problems. Peng \cite{Pe97} introduced the notion of $g$-expectation via solution of backward stochastic differential equation (BSDE for short), which is a dynamic consistent nonlinear expectation. As an application of $g$-expectation, Chen and Epstein \cite{Ce} formulated a continuous-time intertemporal multiple-priors utility with aversion to ambiguity. By Girsanov theorem, it is known that these multiple-priors characterize mean-uncertainty in essence.

Paraller to mean-uncertainty, there is also volatility-uncertainty. In the opinion of Peng \cite{Pe17}, some research on this aspect is equal to defining a manifold on a Euclid space, which is not intrinsic. In fact, Peng \cite{Pe04,Pe05,Pe06} had been trying to establishing the framework of nonlinear expectation space that is just intrinsic to dealing with volatility-uncertainty. Up to now, numerous theoretic results of nonlinear expectation have come to the fore, such as maximal distribution, G-normal distribution, G-Brownian motion, integration with respect to G-Brownian motion, G-martingale, G-martingale representation theorem and G-Girsanov theorem. They are all summarized in Peng's monograph \cite{Pe19}.

Recently, nonlinear expectation theory has been extensively applied in economics, finance and insurance. In fact, Avellaneda et al. \cite{alp95} and  Lyons \cite{Lyo} firstly noted this problem and discussed option pricing in this view. There are still some other explorations along this route. The active issues include asset pricing and stochastic differential utility in \cite{Ej13,Ej14}, risk measure in \cite{pyy21,py21}, option pricing in  \cite{V, Hj18, be19}, interest rate term structure in \cite{ho2101,ho2102,ho2103}, mean–variance hedging in \cite{Bmb}, optimal consumption and portfolio choice in \cite{Lr}, ruin probability in \cite{Sz} and stochastic dominance in \cite{X}, etc.. Based on nonlinear methodology, Gilboa and Schmeidler \cite{GS} not only presented an extension of probability theory, but also obtained a natural generalization of classical expectaion utility theory, i.e., max-min expectation utility theory. Chen and Epstein \cite{Ce} reconsidered the stochastic differential utility model of \cite{GS} with mean-uncertainty. Epstein and Ji \cite{Ej14} reinvestigated the case of \cite{GS} with volatility-ambiguity. In addition, \cite{Ej14} considered the framework of nonlinear expectation put forward by Peng \cite{Pe19}. Due to this point, stochastic differential utility in the framework of nonlinear expectation is called as G-utility, a new notion following G-expectation.

 As a supplement of above applications, a G-utility maximization problem focused on representative investor's terminal wealth is discussed in this paper. Taking a cue from \cite{Wxz}, a similar sufficient condition for optimal strategy is presented. However, the routine of \cite{Wxz} fails since the quadratic variation process of a G-Brownian motion is also a stochastic process. In order to overcome this difficulty, an extended nonlinear expectation space proposed by \cite{Hj14} should be pulled in. In extended nonlinear expectation space, to obtain an optimal strategy of the log-utility optimization problem, a couple of G-FBSDEs related to the optimal strategy is constructed firstly. Then by an auxiliary BSDE, an optimal solution is obtained. A sufficient condition is verified by G-It\^o's formula in the end. The main contributions of this paper are the three pivotal steps of G-martingale approach to solve an optimal strategy in the G-utility maximization problem.

This paper is organized as follows. In Section 2, some preliminaries of nonlinear expectation theory are listed. In Section 3, representative investor's G-utility maximization problem is presented, especially including a sufficient condition for it. In Section 4, an explicit solution of optimal strategy is obtained in the case of log-utility. This special utility function makes the investor's decision to be simplified both in calculation and economic meaning, and also provides direct and important conclusions in the complicated nonlinear setting. In Section 5,  as an application of the main results, a stochastic interest model is considered. Last section gives some concluding remarks. 

\section{Preliminaries}
This section presents an overview of nonlinear expectation theory. It consists of some fundamental conceptions and theorems, such as G-expectation, G-Brownian motion, integration with respect to G-Brownian motion, G-martingale and G-martingale representation theorem and G-Girsanov theorem etc.. The proof of these theorems and further details can be found in \cite{Pe19}. For the sake of simplicity, only the results of one-dimensional case  are recounted in what follows. 

Let $\Omega$ be a given set and $\mathcal{H}$ be a linear space of real valued functions defined on $\Omega$. Assume that $\mathcal{H}$ satisfies the following three conditions,\\

(1) $c\in \mathcal{H}$ for each constant $c \in \mathbb{R}$,\\

(2) $|X| \in \mathcal{H}$, if $ X \in \mathcal{H}$ ,\\

(3) $\phi(X) \in \mathcal{H}$ for any $X \in \mathcal{H}$ and $\phi \in C_{l,lip}(\mathbb{R})$,\\

\noindent where $C_{l,lip}(\mathbb{R})$ denotes the linear space of functions $\varphi$ satisfying local Lipshitz conditions. \\

The definition of sublinear expectation is listed as follows.
\bd\label{ne}
	A sublinear expectation $\mathbb{E}$ is a functional $\mathbb{E}:\mathcal{H}\mapsto \mathbb{R}$ satisfying\\
	
	(1) Monotonicity: 
	$
	\mathbb{E}\left[ X \right] \geqslant \mathbb{E}\left[ Y \right]  ,\ \ \text{~if~} \ X\geqslant Y,\\
	$
	
	(2) Constant preserving: 
	$\mathbb{E}\left[ c \right] =c,\ \ \text{~for~}\ c \in \mathbb{R} ,$\\
	
	(3) Sub-additivity: 
	$\mathbb{E}\left[ X+Y \right]\leqslant \mathbb{E}\left[X \right] +\mathbb{E}\left[ Y \right] ,\ \  X,\ Y \in \mathcal{H}, $\\
	
	(4) Positive homogeneity: 
	$\mathbb{E}\left[ \lambda X \right] =\lambda \mathbb{E} \left[ X \right] ,\ \ \lambda\geqslant 0.$\\
	
\noindent The triple $(\Omega,\mathcal{H},\mathbb{E})$ is called a sublinear expectation space.
\ed

A notable difference between sublinear expectation and linear expectation is (3) in  Definition \ref{ne}. Moreover, their relationship is described by the following representation theorem.

\bt \label{rt}
	Given a nonlinear expectation space $(\Omega,\mathcal{H},\mathbb{E})$, 
there exists a family of linear functionals ${\mathbb E}_{\theta}:\mathcal{H}\mapsto \mathbb{R}$, indexed by $\theta\in \Theta$, such that
	$$\mathbb{E}[X]= \underset{\theta\in \Theta}{\max}{\mathbb E}_{\theta}[X],\ \  X \in \mathcal{H}.$$
\et

From Theorem \ref{rt} , if $\mathbb{E}$ is a sublinear expectation, then the corresponding $E_{\theta}$ is a linear expectation. Accordingly, a conjugate sublinear expectation $\cE$ can be defined by
$$\mathcal{E}[X]=-\mathbb{E}[-X]=\underset{\theta\in \Theta}{\min}{\mathbb E}_{\theta}[X], \ \ X \in \mathcal{H}.$$

In the framework of sublinear expectation space, there are two important distributions, i.e., maximal distribution and G-normal distribution, which are defined as follows. They provide more powerful analysis tools to formulate mean-uncertainty and variance-uncertainty respectively than the classical probabilistic counterparts. 

\bd (1) (maximal distribution) A  random variable $\eta$  on a sublinear expectation space $(\Omega, \mathcal{H}, \mathbb{E})$ is called to be maximally distributed if there exists a bounded, closed and convex subset $\Gamma=[\underline{\mu},\overline{\mu}] \subset \mathbb{R}$ such that
	$$
	\mathbb{E}[\varphi(\eta)]=\max _{y \in \Gamma} \varphi(y),\ \  \varphi \in C_{l, L i p}\left(\mathbb{R}\right).
	$$
 Particularly, the maximal distribution on $[\underline{\mu},\overline{\mu}]$ is denoted by $N([\underline{\mu},\overline{\mu}],0)$.
	
	(2) (G-normal distribution) A  random variable $X$ on a sublinear expectation space $(\Omega, \mathcal{H}, \mathbb{E})$ is called to be $G$-normally distributed if for each $a, b \geq 0$, \\
	$$
	a X+b \bar{X} \stackrel{d}{=} \sqrt{a^{2}+b^{2}} X
	$$
	where $\bar{X}$ is an independent copy of $X$, i.e., $\bar{X} \stackrel{d}{=} X$ and $\bar{X} \perp X$. Particularly, the zero mean G-normal distribution on  $[\underline{\sigma}^2,\overline{\sigma}^2]$ is denoted by $N(0,[\underline{\sigma}^2,\overline{\sigma}^2])$.
\ed

  Peng \cite[Definition 3.1.2]{Pe19} defined a G-Brownian motion.

\bd(G-Brownian motion)
	A  stochastic process $(B_{t})$ on a sublinear expectation space $(\Omega, \mathcal{H}, \mathbb{E})$ is called a $G$-Brownian motion if \\
	
	(i) $B_{0}(\omega)=0$,\\
	
	(ii) For each $t, s \in [0,T], \ B_{t+s}-B_{t}$ and $B_{s}$ are identically distributed,\ and $B_{t+s}-B_{t}$ is independent of  $\left(B_{t_{1}}, B_{t_{2}}, \cdots, B_{t_{n}}\right)$,\  for each $n \in \mathbb{N}$ and $0 \leq t_{1} \leq \cdots \leq t_{n} \leq t$,\\
	
	(iii) $\lim _{t \downarrow 0} \mathbb{E}[\left|B_{t}\right|^{3}] t^{-1}=0$.\\
	
Furthermore, if $\mathbb{E}\left[B_{t}\right]=\mathbb{E}\left[-B_{t}\right]=0$, then $(B_{t})$ is called a symmetric $G$-Brownian motion.
\ed

 Let $\Omega_{T}=C_{0} ([0, T];\mathbb{R} )$ denote the space consisting of real-valued continuous functions with initial values equal to zero. It is equipped with the supreme norm on $[0, T]$. Denote  
$$Lip(\Omega_T)=\left\{\varphi\left(B_{t_{1}}, B_{t_{2}}, \cdots, B_{t_{n}}\right),\ 0=t_0\leq t_1\leq t_2\leq \cdots \leq t_n=T, \ \varphi \in C_{l,lip}(\mathbb{R})\right\}.$$
Peng \cite[Page 53-54]{Pe19} established G-expectation and conditional G-expectation respectively on $(\Omega, Lip(\Omega_T))$. The former is a sublinear expectation such that the canonical process $(B_t)$ is a G-Brownian motion. The latter is defined as follows, 

$$
\begin{array}{l}
	\mathbb{E}_{t}[\varphi(B_{t_{1}}-B_{t_{0}}, B_{t_{2}}-B_{t_{1}}, \cdots, B_{t_{m}}-B_{t_{m-1}})]=\tilde{\varphi}(B_{t_{1}}-B_{t_{0}}, B_{t_{2}}-B_{t_{1}}, \cdots, B_{t_{i}}-B_{t_{i-1}}),
\end{array}
$$
where
$
\tilde{\varphi}\left(x_{1}, \cdots, x_{i}\right)=\mathbb{E}\left[\varphi\left(x_{1}, \cdots, x_{i}, B_{t_{i+1}}-B_{t_{i}}, \cdots, B_{t_{m}}-B_{t_{m-1}}\right)\right]
,\ \ t=t_i~(1\leq i\leq m, \ m\leq n)$.

Then It\^o integration with respect to G-Brownian motion was defined similar to the case with respect to Brownian motion. Distinguished from the classical one, the quadratic variation process of a G-Brownian motion, i.e., 
$$ \langle B\rangle_{t} = B_{t}^{2}-2 \int_{0}^{t} B_{s} \dif B_{s}, \ \ t \in [0,T] $$ 
is an increasing process with $\langle B\rangle_{0}=0$. Furthermore, $(\langle B\rangle_{t})$ characterizes some statistic uncertainty of the $G$-Brownian motion  $(B_{t})$ . For the degenerated case $\underline{\sigma}^2=\overline{\sigma}^2$ and any $t$, $\langle B\rangle_{t}$ always obeys maximal distribution  $N([\underline{\sigma}^2 t, \overline{\sigma}^2 t],0)$. 

Given a G-Brownian motion $(B_t)$ and its quadratic variation process $(\langle B_t \rangle)$ , Peng \cite[Page 71]{Pe19} proved G-It\^o's formula for the following G-It\^o process $(X_t)$, 
$$X_t=X_0+\int_0^t \alpha_s \dif s+\int_0^t \beta_s\dif B_s+\int_0^t \eta_s \dif \langle B \rangle_{s},~ t\in[0,T],  $$
where $(\alpha_t) \ \text{and} \  (\eta_t)$ belong to $ M_G^1(0,T)$ which denotes the class of integrable stochastic processes in the sense of G-expectation. And $(\beta_t) $ is in $M_G^2(0,T)$, which collects exclusively square integrable stochastic processes.
\bt
	Let $\Phi$ be a $C^{2}$-function on $\mathbb{R}$ such that $\partial_{x x}^{2} \Phi$ satisfies the polynomial growth condition, then for each $t \in [0,T], \ X_{t} \in L_{G}^{2}(\Omega_{t})$, the following formula holds,
	$$
	\begin{aligned}
		\Phi\left(X_{t}\right)-\Phi\left(X_{0}\right)&=\int_0^t \partial_x \Phi(X_u)\alpha_u \dif u+\int_0^t \partial_x \Phi(X_u)\beta_u \dif B_u\\
		&+\int_0^t \left[\partial_x \Phi(X_u)\eta_u+\frac{1}{2}\partial_{xx}^2 \Phi(X_u)\beta_u^2 \right] \dif \langle B \rangle_u,
	\end{aligned}
	$$
	where  $L_{G}^{2}\left(\Omega_{t}\right)$ is the completion of the space $Lip(\Omega_t)$ under $L^2$-norm of G-expectation.
\et

To present the G-Girsanov theorem, the notion of  G-martingales and  related martingale representation theorem are listed previously. 
\bd(G-martingale)
	A process $(M_{t})$ is called a G-martingale if for any $t \in [0,T],\  M_{t} \in L_{G}^{1}(\Omega_{t})$ and for any $s \in[0, t]$, it holds
	$$
	\mathbb{E}[M_{t} \mid \Omega_{s}] = M_{s},\ \ q.s..
	$$
	Furthermore, if $M$ also satisfies
	$$
	\mathbb{E}[-M_{t} \mid \Omega_{s}]=-M_{s},\ \ q.s.,
	$$
	then it is called a symmetric G-martingale.
\ed

Here, $q.s.$ is the abbreviation for quasi-surely, similar to $a.s.$ in stochastic analysis.

\bt(Representation theorem of G-martingale)\label{gmrp}
	Let $p>1$. For some suitable $\xi \in \mathbb{L}_{G}^{p}$ and denoting $Y_{t}:=\mathbb{E}_{t}^{G}[\xi]$, there exist unique $Z \in \mathcal{H}_{G}^{p}$ and $\eta \in \mathbb{M}_{G}^{*}$ satisfying
	\begin{equation*}
		Y_{t}=Y_{0}+\int_{0}^{t} Z_{s} \dif  B_{s}-\left[\int_{0}^{t} G\left(\eta_{s}\right) \dif s-\frac{1}{2} \int_{0}^{t} \eta_{s} \dif \langle B\rangle_{s}\right],\ \ q.s.,\ \ t \in [0,T].
	\end{equation*}  

\et
 The further details of the theorem can be retrieved in \cite{Pe14}. 

To claim G-Girsanov theorem, we introduce an auxiliary extended $\tilde{G}$-expectation space $(\tilde{\Omega}_{T}, L_{\tilde{G}}^{1}(\tilde{\Omega}_{T}), \mathbb{E}^{\tilde{G}})$ with $\tilde{\Omega}_{T}=C_{0}([0, T], \mathbb{R}^{2})$ and
$$
\tilde{G}(A):=\frac{1}{2} \sup _{\underline{\sigma}^{2} \leq v \leq \bar{\sigma}^{2}} \operatorname{Tr}\left[A\left[\begin{array}{cc}
	v & 1 \\
	1 & v^{-1}
\end{array}\right]\right], \quad A \in \mathbb{S}_{2},
$$
where $\mathbb{S}_{2}$ represents the class of $2\times 2$ symmetric matrices.
Moreover, it is clear that the canonical process  $(B_{t}, \tilde{B}_{t})$  is a 2-dimensional G-Brownian motion in the extended space $(\tilde{\Omega}_{T}, L_{\tilde{G}}^{1}(\tilde{\Omega}_{T}), \mathbb{E}^{\tilde{G}})$.

\bt \label{GT}(G-Girsanov Theorem, \cite{Hj14})
	Let $(b_{t})$ and $(d_{t})$ be bounded and adapted processes. Then $\bar{B}_{t}:= B_{t}-\int_{0}^{t} b_{s} \dif s - \int_{0}^{t} d_{s} \dif \langle B \rangle_{s} $ is a G-Brownian motion under $\tilde{\mathbb{E}}$,
	where $\tilde{\mathbb{E}}$ satisfies
	\begin{equation*}
		\begin{aligned}
			\tilde{\mathbb{E}}_{t}[\cdot] = \mathbb{E}_{t}^{\tilde{G}} & {\left[~\cdot \exp \left \lbrace \int_{t}^{T} d_{s} \dif B_{s}-\dfrac{1}{2} \int_{t}^{T} |d_{s} |^{2} \dif \langle B \rangle_{s} - \int_{t}^{T} b_{s} d_{s} \dif s \right. \right.} \\
			&\left. \left. ~ + \int_{t}^{T} b_{s} \dif \tilde{B}_{s} - \dfrac{1}{2} \int_{t}^{T} |b_{s}|^{2} \dif \langle \tilde{B} \rangle_{s} \right \rbrace \right],\ \ t \in [0,T].
		\end{aligned}
	\end{equation*}
\et
\br\label{228} There  exists a stochastic process $(\rho_t) \in M_{\tilde{G}}^{2,+} (\tilde{\Omega}_{T}) $, which denotes a set of the whole nonnegative ones in $M_{\tilde{G}}^2(\tilde{\Omega}_{T})$, such that $d \tilde{B}_t=\rho_t \dif B_t,\ t\in[0,T],$ holding for some G-Brownian motions because of the special G-covariance of $(B_t)$ and $(\tilde{B}_t)$. In fact, assume that a process $(B_t^{\prime})$ satisfies
	\begin{equation*}
		B_t^{\prime}=\int_0^t \gamma_s \dif W_s, \ \  t\in[0,T],
	\end{equation*}
where $(W_t)$ is a standard Brownian motion  and  $\gamma_t$ is maximally distributed as $N([\underline{\sigma}^2 t,\overline{\sigma}^2 t],0)$ for any $t\in[0,T]$. It is clear that $(B_t^{\prime})$ is a G-Brownian motion. Meanwhile, if $(\tilde{B}_t^{\prime})$ is defined by $ \displaystyle \tilde{B}_t^{\prime}= \int_0^t \dfrac{1}{\gamma_s} \dif W_s,\ t\in[0,T]$, then there exists  $\rho_t = \dfrac{1}{{\gamma_t}^2} \in L_{\tilde{G}}^{2,+}(\tilde{\Omega}_{T}) $ such that  $d \tilde{B}_t=\rho_t \dif B_t$ holding.

To obtain explicit solution of an optimal investment strategy, without loss of generality, a class of 2-dimensional G-Brownian motions $\{(B_t^{\prime},\tilde{B}_t^{\prime})\}$  are just considered in this paper .
\er

Under suitable conditions, Peng \cite[Page 103,105]{Pe19} proved the existence and uniqueness of G-SDEs for $(X_t)$ and G-BSDEs for $(Y_t)$ which are both driven by G-Brownian motion, i.e., 
\begin{equation*}
	X_{t}=X_{0}+\int_{0}^{t} b\left(s, X_{s}\right) \dif s+\int_{0}^{t} h\left(s, X_{s}\right) \dif\langle B\rangle_{s}+\int_{0}^{t} \sigma\left(s, X_{s}\right) \dif B_{s},  \ \  t \in[0, T],
\end{equation*}
and
\begin{equation*}
	Y_{t}=\mathbb{E}\left[\xi+\int_{t}^{T} f\left(s, Y_{s}\right) \dif s+\int_{t}^{T} h\left(s, Y_{s}\right) \dif \langle B\rangle_{s} \mid \Omega_{t}\right], \ \  t \in[0, T].
\end{equation*}
The couple of a G-SDE and a G-BSDE are called G-FBSDEs.

\section{Problem Formulation of Investor's G-Utility Maximization}

Consider a continuous-time financial market with fixed time horizon $[0,T]$. There are a bond and a stock, regarded as safe and risk asset respectively. Representative investor aims to maximize expected utility with respect to portfolio wealth at the terminal date $T$. As usual, the frictionless conditions of financial market are satisfied. For instance, continuous trade and negative transaction are allowed; there are no transaction cost or tax; and all assets are divisible. Uncertainty directly affects the investor's views on the price behavior of risky asset. Significantly, his own expectation utility and optimal decision behavior would be influenced in an unusual way.

The uncertainty source from financial market is characterized by the price process  $(S_t)$ of the stock driven by the G-Brownian Motion $(B_t)$. Together with the price process $(R_t)$ of the bond, both are specified as
\begin{align*}
	\dif R_t &= rR_t \dif t, \ \ \ \  R_0=1,\\
	\dif S_t &= S_t\big(\mu \dif t + \sigma \dif B_t + c \dif \langle B \rangle_t \big), \ \ \ \  S_0=1,
\end{align*}
where risk-free rate $r>0$, appreciation rate $\mu>r$, volatility $\sigma >0$, and a constant $c>0 $. Due to the appearance of the quadratic variation process $( \langle B \rangle_t )$, interval $ [ c\underline{\sigma}^2, c\overline{\sigma}^2 ]$ has a special meaning that indicates the mean-uncertainty implicit in the appreciation rate of the risk asset. 

The next concepts are associated with investor's portfolio decision. Let $\pi=(\pi_{t})$ be a trading strategy, where $\pi_t$ is the proportion invested in stock at any $t \in [0,T]$, and only depending on the information before time $t$.     
\bd
	A trading strategy $\pi$ is called to be admissible if $ \pi \in M_{\tilde{G}}^2 (\Omega_T). $
\ed
Denote $\Pi$ the collection of all admissible trading strategies. For an element $\pi \in \Pi$ and initial capital $x_0$, wealth process $X^{x_0,\pi}$ of the investor follows
\begin{equation*}
	\begin{cases}
		\dif X_t^{x_0,\pi} = X_t^{x_0,\pi}r\dif t + \pi_t X_t^{x_0,\pi}\left((\mu-r)\dif t + \sigma \dif B_t+c \dif \langle B \rangle_t\right),\ \  t \in [0,T],\\
		X_0^{x_0,\pi}=x_0.
	\end{cases}
\end{equation*}
For convenience to calculate, a discounted wealth process is considered and defined by $\bar{X}_t^{x_0,\pi} := X_t^{x_0,\pi}{R_t}^{-1}$. By G-It\^o's formula, one has
\begin{equation*}\label{bar_X}
	\begin{cases}
		\dif \bar{X}_t^{x_0,\pi}=\pi_t \bar{X}_t^{x_0,\pi} ((\mu-r)\dif t +\sigma \dif B_t+c \dif \langle B\rangle_t),\ \  t \in [0,T],\\
		\bar{X}_0^{x_0,\pi}=x_0.
	\end{cases}
\end{equation*}
That is,
\begin{equation*}
	\bar{X}_t^{x_0,\pi}=x_0\exp\left\{\int_0^t \pi_s\left(\mu-r\right)\dif s+\int_0^t \pi_s \sigma \dif B_ s+\int_0^t \left(\pi_s c-\frac{1}{2}\pi_s^2 \sigma^2\right)\dif \langle B \rangle_s\right\},\ \  t \in [0,T].
\end{equation*}

As a price-taker, the investor's decision behaves uncertainty. There are varieties of characterization involved in objective function specification of optimization procedures. It includes both of the selection in uncertainty measure and utility function of all potential outcomes. Sublinear expectation $\mE$ or its conjugate measure seems to be a suitable choice, but encounters a serious challenge for solving the optimization model. The key step is to adopt the auxiliary extended sublinear expectation space $\mE^{\tilde{G}}$ introduced in Section 2. This technique makes it possible to exploit G-martingale approach to solve G-FBSDEs associated with the optimization problem. For characterizing uncertainty-aversion  of the investor, the conjugate measure $\cE$  is a more ideal one, defined as
$$\cE [\cdot]:=-\mE^{\tilde{G}} [-~\cdot].$$ 

By the definition of $\mE^{\tilde{G}}$, the below properties hold,\\

(1) $\cE[c]=-\mE^{\tilde{G}}[-c]=c$,\ \ for all $c\in \mR$,\\

(2) $\cE[X]=-\mE^{\tilde{G}}[-X]\geq -\mE^{\tilde{G}} [-Y]=\cE [Y]$,\ \ for all $X, Y \in \cH$ and $X \geq Y$,\\

(3) $\cE[\lambda X]=-\mE^{\tilde{G}} [-\lambda X]=\lambda \left(-\mE^{\tilde{G}} [-X]\right)=\lambda \cE [X]$,\ \ for all $\lambda \in \mR^+$,\\

(4) $\cE[X]-\cE[Y]=-\left(\mE^{\tilde{G}} [-X]-\mE^{\tilde{G}} [-Y]\right)\geq -\mE^{\tilde{G}} [-(X-Y)]=\cE[X-Y]$,\ \ for all $X,Y\in \cH$.\\

Assume that the investor has a utility function $U(~\cdot)$ of the terminal discounted wealth. Then a desired objective function will be specified straightly. It is called G-utility as in the following maximal aim,
\begin{equation}\label{op}
	\max_{\pi \in \Pi}\cE \left[U\left(\bar{X}_T^{x_0,\pi}\right)\right].
\end{equation}
For two strategies $\pi^* \in \Pi$ and $ \pi \in \Pi$ , by the concavity of $U$, the following inequalities hold,
\begin{equation}\label{sc1}
	\begin{aligned}
		&\cE[U(\bar{X}_T^{x_0,\pi^*})]-\cE[U(\bar{X}_T^{x_0,\pi})]\\
		&\geq \cE [U(\bar{X}_T^{x_0,\pi^*})-U(\bar{X}_T^{x_0,\pi})]\\
		&=\cE [U(\bar{X}_T^{x_0,\pi^*})-U(\bar{X}_T^{x_0,\pi^*})-U^{\prime}(\bar{X}_T^{x_0,\pi^*})(\bar{X}_T^{x_0,\pi}-\bar{X}_T^{x_0,\pi^*})-U^{\prime \prime}(\bar{X}_T^{x_0,\pi^*})(\bar{X}_T^{x_0,\pi}-\bar{X}_T^{x_0,\pi^*})^2]\\
		&\geq \cE [U^{\prime}(\bar{X}_T^{x_0,\pi^*})(\bar{X}_T^{x_0,\pi^*}-\bar{X}_T^{x_0,\pi})].
	\end{aligned}
\end{equation}

In the end of this section, a sufficient condition of G-utility maximization problem (\ref{op}) is given. It is a prerequisite to apply G-martingale approach. According to (\ref{sc1}), one can verify this point directly.
\bp\label{sc} 
	If there exists a trade strategy $\pi^* \in \Pi$, such that
	\begin{equation}\label{condition}
		U^{\prime}(\bar{X}_T^{x_0,\pi^*})\bar{X}_T^{x_0,\pi} \text {is a symmetric G-martingale for any } \pi \in \Pi,
	\end{equation} 
	then $\pi^*$ is optimal.
\ep
Note that symmetric G-martingales play a very critical role in asset pricing in nonlinear expectation theory. Referring to Beissner \cite{be19}, it is a fair valuation measure and its existence is equivalent to no-arbitrage condition under uncertainty.
 
\section{an Explicit Solution with Log-Utility}

Given condition (\ref{condition}) and log-utility function, i.e., $U(x)=\ln x$, an explicit solution of optimal strategies can be derived. By Proposition \ref{sc}, an optimal strategy $\pi^*$ is required to satisfy that
\begin{equation}\label{ecmv}
	\begin{aligned}
		&\dfrac{1}{\bar{X}_T^{x_0,\pi^*}} x_0 \exp \left\{\int_0^t \pi_s (\mu - r) \dif s + \int_0^t \pi_s \sigma \dif B_t + \int_0^t \Big( \pi_s c - \dfrac{1}{2}\pi_s^2 \sigma^2 \Big) \dif  \langle B \rangle_t \right\}\\
		&\text{is a symmetric martingale for any }  \pi \in \Pi. 
	\end{aligned} 
\end{equation}

Here, the optimal discounted terminal wealth $\bar{X}_T^{x_0,\pi^*}$ represents a baseline result of rational decision in financial market. The ratio of $x_0$ to $\bar{X}_T^{x_0,\pi^*}$ has the special meaning of discount factor to some extent just like the classical case in which $(S_t)$ is driven by Brownian motions. Hence, set 
\begin{equation}\label{Z}
	Z_t^*=\frac{\cE_t\left[\dfrac{1}{\bar{X}_T^{x_0,\pi^*}}\right]}{\cE\left[\dfrac{1}{\bar{X}_T^{x_0,\pi^*}}\right]}, \ \ \  t\in [0,T].
\end{equation}
Then 
\begin{equation*}
	Z_T^*=\dfrac{1}{\cE \left[ \dfrac{1}{\bar{X}_T^{x_0,\pi^*}} \right] \bar{X}_T^{x_0,\pi^*}}
	, \  \  \ Z_{\tau}^* = \dfrac{\cE_{\tau} \left[ \dfrac{1}{\bar{X}_T^{x_0,\pi^*}} \right] }{\cE \left[ \dfrac{1}{\bar{X}_T^{x_0,\pi^*}} \right] },
\end{equation*}
for any stopping time $\tau \leq T$,\  q.s.. A further result is summarized in the following lemma.
\bl \label{G_SDE}
	If $\pi^*\in \Pi $ satisfies condition (\ref{ecmv}), then the triplet $(\bar{X}^{x_0,\pi^*},\bar{\pi}^*,\bar{Z}^*)$ solves  the following linear G-FBSDEs,
	\begin{equation}\label{fbsde}
		\begin{cases}
			\dif \bar{X}_t^{x_0,\pi} = \pi_t\bar{X}_t^{x_0,\pi} \big( (\mu - r)\dif t + \sigma \dif B_t + c \dif \langle B \rangle_t \big),\ \ t\in [0,T],\\
			X_0 = x_0,\\
			\dif Z_t = -\dfrac{c}{\sigma}Z_t \dif B_t - \dfrac{\mu - r}{\sigma} Z_t \dif \tilde{B}_t,\ \ t\in [0,T],\\
			Z_T^* = \dfrac{1}{\cE \left[ \dfrac{1}{\bar{X}_T^{x_0,\pi^*}} \right] \bar{X}_T^{x_0,\pi^*}}.
		\end{cases}
	\end{equation}
\el
\begin{proof}
	Suppose $\pi^*$ satisfies condition (\ref{ecmv}). It is clear that $\bar{X}^{x_0,\pi^*}$ solves the G-SDE for $\bar{X}$ and $Z^*$ is a symmetric G-martingale. For any stopping time $\tau \leq T$, set $\pi_{t}^{\tau} = I_{\{t\leq \tau\}}$, then $\pi^{\tau} \in \Pi$.
	
	Substituting $\pi^{\tau}$ into (\ref{ecmv}), one obtains
	\begin{align*}
		\hspace{1em} \cE\left[\dfrac{1}{\bar{X}_T^{x_0,\pi^*}}\right] & \cE\left[Z_T^*x_0 \exp \left\{\int_0^T \pi_s\left(\mu - r\right)\dif s + \int_0^T \pi_s \sigma \dif B_s + \int_0^T \Big(\pi_s c - \dfrac{1}{2}\pi_s^2 \sigma^2 \Big) \dif \langle B \rangle_s \right\} \right] \\
		&= \cE\left[\dfrac{1}{\bar{X}_T^{x_0,\pi^*}}\right]\cE\left[\cE_{\tau}[Z_T^*]x_0 \exp\left\{\int_0^{\tau} \!\!(\mu-r)\dif s + \! \int_0^{\tau} \! \sigma \dif B_s + \! \int_0^{\tau} \Big( c-\frac{1}{2} \sigma^2\Big)\dif \langle B \rangle_s\right\}\right]\\
		&= \cE\left[\dfrac{1}{\bar{X}_T^{x_0,\pi^*}}\right]\cE\left[Z_{\tau}^* x_0 \exp\left\{\int_0^{\tau} \left(\mu - r\right)\dif s + \int_0^{\tau} \sigma \dif B_s + \int_0^{\tau} \Big( c - \dfrac{1}{2} \sigma^2 \Big)\dif \langle B \rangle_s\right\}\right].
	\end{align*}

\noindent Combining with condition ($\ref{ecmv}$), it is clear that
	\begin{equation*}\label{ecmv1}
		Z_{t}^* x_0\exp \left \{\int_0^{t} \left(\mu - r\right)\dif s + \int_0^{t}  \sigma \dif B_s + \int_0^{t} \Big( c - \dfrac{1}{2} \sigma^2\Big) \dif \langle B \rangle_s \right \} \text{is a symmetric G-martingale}.
	\end{equation*}
	By (\ref{Z}), one can claim that $Z_t^*$ is a G-martingale. With Theorem  \ref{GT}, one derives
	\begin{equation}\label{zr}
	\dif Z^*_t = -\dfrac{c}{\sigma}Z^*_t \dif B_t - \dfrac{\mu-r}{\sigma} Z^*_t \dif \tilde{B}_t, \ \ t\in[0,T].
	\end{equation}
	Meanwhile, it is direct to verify  $Z_T^*=\dfrac{1}{\cE\left[\dfrac{1}{\bar{X}_T^{x_0,\pi^*}}\right]\bar{X}_T^{x_0,\pi^*}}$ .
\end{proof} 

Lemma \ref{G_SDE} states that $\pi^*$ is a solution of (\ref{fbsde}). If there exists a solution to (\ref{fbsde}) satisfying condition (\ref{condition}), then it is an optimal strategy indeed. This presumption matches perfectly with optimization problem (\ref{op}) adopting log-utility. The result is presented in the next theorem. 

\bt\label{t_1}
	For $ U(x) = \ln x $, \  $\pi^* = (\pi_t^*) = \Big( \dfrac{(\mu-r)\rho_t + c}{\sigma^2} \Big) $ is an optimal strategy. 
\et
\begin{proof}
	Define an auxiliary process $(M_t)$ as follows
	\begin{equation}\label{A1}
		\begin{cases}
		  ~ \dif M_t=a_t M_t \dif t+b_t M_t \dif \langle B\rangle_t, \ \ t\in [0,T],\\
		  ~ M_T=1.
		\end{cases}
	\end{equation}
By G-It\^o's formula,
	\begin{equation}\label{solve}
		\begin{aligned}
			\dif \Big ( M_t \dfrac{1}{\bar{X}_t^{x_0, \pi^*}} \Big ) =& ~ M_t \dfrac{1}{\bar{X}_t^{x_0, \pi^*}} \big[a_t - \pi_t^*(\mu - r)\big] \dif t + M_t \dfrac{1}{\bar{X}_t^{x_0, \pi^*}} \left[ b_t + \pi_t^*(\pi_t^* \sigma^2 -c) \right] \dif \langle B \rangle_t\\
			& - M_t \dfrac{\pi_t^* \sigma}{\bar{X}_t^{x_0,\pi^*}}\dif B_t,\ \ t\in [0,T].
		\end{aligned}
	\end{equation}
To ensure that $\Big( M_t \dfrac{1}{\bar{X}_t^{x_0,\pi^*}}\Big)$ is a symmetric G-martingale, it suffices to set
	\begin{equation}\label{a,b}
		\begin{cases} 
			~ a_t=\pi^*_t(\mu-r),\\ ~ b_t=-\pi_t^*(\pi_t^* \sigma^2-c).
		\end{cases}
	\end{equation}
Substituting (\ref{a,b}) into (\ref{solve}) and combining with (\ref{A1}), one derives
	\begin{equation}\label{mx}
		M_t \dfrac{1}{\bar{X}_t^{x_0,\pi^*}}=\dfrac{1}{\bar{X}_T^{x_0,\pi^*}}+\int_t^T M_s \dfrac{\pi_s^* \sigma}{\bar{X}_s^{x_0,\pi^*}}\dif B_s, \ \ t \in [0,T].
	\end{equation}

Note that $ \displaystyle \Big(\int_0^t M_s \dfrac{\pi_s^* \sigma}{\bar{X}_s^{x_0,\pi^*}}\dif B_s \Big ) $ is a symmetric G-martingale, by taking conditional expectation with respect to two sides of equation, one obtains
	\begin{equation*} 
		M_t \dfrac{1}{\bar{X}_t^{x_0,\pi^*}}=\cE_t\left[\dfrac{1}{\bar{X}_T^{x_0,\pi^*}}\right], \ \ t \in [0,T].
	\end{equation*}
Together with (\ref{mx}), this implies
	\begin{equation*}
		\dif\cE_t\left[\dfrac{1}{\bar{X}_T^{x_0,\pi^*}}\right] = \dif \Big (M_t \dfrac{1}{\bar{X}_t^{x_0,\pi^*}} \Big ) = - \cE_t\left[\dfrac{1}{\bar{X}_T^{x_0,\pi^*}}\right]\pi_t^* \sigma \dif B_t, \ \ t \in [0,T].
	\end{equation*}
Combining with (\ref{Z}), it is clear that
	\begin{equation*}
		\dif Z_t^* = -Z_t^*\pi_t^* \sigma \dif B_t, \ \ t \in [0,T].
	\end{equation*}
By Theorem \ref{gmrp}, Remark \ref{228} and (\ref{zr}), one concludes that
	\begin{equation}\label{ops}
		\pi^*_t=\frac{(\mu-r)\rho_t+c}{\sigma^2}, \ \ t \in [0,T].
	\end{equation}
	
Finally, it can be proved that the solution (\ref{ops}) truly satisfies condition ($\ref{condition}$). For any another strategy $\pi$, consider the ratio of $ \bar{X}_t^{x_0,\pi} $ to $ \bar{X}_t^{x_0,\pi^*} $ by G-It\^o's formula. 

	\begin{equation*}
		\begin{aligned}
			\dif \dfrac{\bar{X}_t^{x_0,\pi}}{\bar{X}_t^{x_0,\pi^*}} =& ~ \dfrac{\dif\bar{X}_t^{x_0,\pi}}{\bar{X}_t^{x_0,\pi^*}} + \bar{X}_t^{x_0,\pi} \dif \dfrac{1}{\bar{X}_t^{x_0,\pi^*}} - \dfrac{c \pi_t \bar{X}_t^{x_0,\pi}}{\bar{X}_t^{x_0,\pi^*}} \dif \langle B \rangle_t - \dfrac{ (\mu-r) \pi_t \bar{X}_t^{x_0,\pi} }{\bar{X}_t^{x_0,\pi^*}} \dif t \\
			=& ~ \dfrac{\pi_t\bar{X}_t^{x_0,\pi}}{\bar{X}_t^{x_0,\pi^*}} \Big[ (\mu-r)\dif t + c \dif \langle B\rangle_t - (\mu-r)\dif t - c \dif\langle B\rangle_t \Big] \\
			&+ \dfrac{\pi^*_t\bar{X}_t^{x_0,\pi}}{\bar{X}_t^{x_0,\pi^*}} \Big[(\mu-r)\dif t - (\mu-r) \dif t \Big] + \dfrac{(\sigma \pi_t - c/\sigma)\bar{X}_t^{x_0,\pi}}{\bar{X}_t^{x_0,\pi^*}} \dif B_t - \dfrac{(\mu-r)\bar{X}_t^{x_0,\pi}}{\sigma \bar{X}_t^{x_0,\pi^*}} \dif \tilde{B}_t \\
			=& ~\dfrac{(\sigma \pi_t - c/\sigma)\bar{X}_t^{x_0,\pi}}{\bar{X}_t^{x_0,\pi^*}} \dif B_t - \dfrac{(\mu-r)\bar{X}_t^{x_0,\pi}}{\sigma \bar{X}_t^{x_0,\pi^*}} \dif \tilde{B}_t, \ \ t\in[0,T].
		\end{aligned}
	\end{equation*}
So $ \Big (\dfrac{\bar{X}_t^{x_0,\pi}}{\bar{X}_t^{x_0,\pi^*}} \Big )$ is a symmetric G-martingale. And $(\pi^*_t)$ in (\ref{ops}) ensures condition (\ref{condition}). It also means that $\pi^*$ is an optimal strategy.
\end{proof}

Up to now, the G-utility maximization model with log-utility is solved explicitly by G-martingale approach. As far as the optimal solution (\ref{ops}) is concerned, volatility uncertainty appears to be absorbed into $\rho_t$,  which obeys maximal distribution $N\big([{\underline \rho}t,{\overline \rho}t],0\big)$ according to Remark \ref{228}. So, some economic explanations are obtained straightly. Firstly, for usual condition $\mu \geq r$, one know that the optimal strategy  increases with the increasing of $\underline \rho$ or $\overline \rho $ . Secondly, the same increasing relationship holds with respect to $c$, since the parameter characterizes the degree of mean increases. Thirdly, the optimal strategy is decreasing with the increasing of $\sigma$, which means that the investors are averse to the degree of volatility uncertainty. All of the three explanations are consistent with our intuitions.

In the end of this section, two generalized results are given. One is a time-varying coefficient model and the other is a model with dividend. Because of the proof similar to Theorem \ref{t_1}, it is omitted.
\bc
	(1) Consider that stock price process follows$$\dif S_t= S_t\big(\mu_t \dif t+ \sigma_t \dif B_t + c_t \dif\langle B \rangle_t\big),\ \ t\in[0,T],$$ where $\mu,c \text{~and~} \sigma$ are all uniformly bounded and adapted processes. And assume that there exists $\varepsilon >0$ such that $\sigma_t \geq \varepsilon$. Then $\pi^* = \Big ( \dfrac{(\mu_t - r)\rho_t + c_t}{\sigma_t^2} \Big ) $ is an optimal strategy in log-utility case. \\
	
	(2) Consider that stock price process follows
	$$\dif S_t= S_t \Big[\Big(\mu+\sum_{i=0}^N d_i\delta_{t_i}(t)\Big) \dif t+ \sigma \dif B_t + c \dif \langle B \rangle_t \Big],\ \ t\in [0,T],$$
	where $d_i$ and $\delta_{t_i}$ denote the dividend and Dirac function respectively at time $t_i \in [0,T] $. Then $\pi^* = \Big ( \dfrac{\big (\mu + d_i\delta_i(t) - r \big )\rho_t + c}{\sigma^2}\Big ) $ is an optimal strategy in the case of log-utility.
\ec

\section{Generalization to the stochastic interest model}

This section will go straightly to discuss how random interest rate affects optimal strategy. It means there exists no real safety in financial market just as in reality. Return of safe asset also behaves uncertainty. Here the stochastic interest rate model proposed by Dothan \cite{Do} is considered and rewritten as 
\begin{equation*}
	\dif \tilde{R}_t = \tilde{R}_t (r\dif t + \sigma_r \dif B_t),\ \  t\in[0,T],
\end{equation*}
where the appreciation rate $r>0$,  volatility of bond is denoted by $\sigma_r$. And $\left(B_t\right)$ is G-Brownian motion. Recall that both the stock price process $(S_t)$ and the admissible strategy set $\Pi$ are defined as before.
Similarly, for an admissible trading strategy $\pi$ and an initial capital $x_0$, the wealth process $\tilde{X}^{x_0,\pi}$ is characterized as 
\begin{equation*}
	\begin{cases}
		~\dif\tilde{X}_t^{x_0,\pi}=\tilde{X}_t^{x_0,\pi}\left(r\dif t+\sigma_r \dif B_t\right)+\pi_t \tilde{X}_t^{x_0,\pi}((\mu-r)\dif t +\left(\sigma-\sigma_r\right) \dif B_t+c \dif \langle B\rangle_t),~t\in[0,T],\\
		~\tilde{X}_0^{x_0,\pi}=x_0.
	\end{cases}
\end{equation*}
And its solution is
\begin{equation}\label{X_tilde}
	\begin{aligned}
	\tilde{X}^{x_0,\pi}_t = & x_0 D_t \exp \left \{ \int_0^t \pi_s(\mu - r) \dif s + \int_0^t \pi_s(\sigma - \sigma_r) \dif B_s ~- \right.\\
	& \left. \int_0^t \Big [\frac{1}{2}\pi_s^2 (\sigma - \sigma_r )^2 + \pi_s\sigma_r (\sigma - \sigma_r) - c\pi_s \Big ] \dif \langle B\rangle_s \right \},\ \ t\in[0,T].
	\end{aligned}
\end{equation}
where $ D_t : = \exp \Big \{ rt + \sigma_r B_t - \dfrac{1}{2}\sigma_r^2 \langle B \rangle_t \Big \} ,\ \ t\in[0,T] $.

In the case of $U(x)=\ln x$, the investor's objective becomes

\begin{equation}\label{saim}
	\max_{\pi \in \Pi}\cE [\ln (\tilde{X}_T^{x_0,\pi})].
\end{equation}
Similar to previous arguments, the following sufficient condition of optimal strategy can be derived.
\bp\label{sec}
	If there is a trade strategy $\pi^* \in \Pi$ such that
	\begin{equation}\label{scondition}
		U^{\prime}(\tilde{X}_T^{x_0,\pi^*})\tilde{X}_T^{x_0,\pi} \text {is a symmetric G-martingale for any } \pi \in \Pi,
	\end{equation} 
	then $\pi^*$  is the optimal one.
\ep
To get optimal strategy $\pi^*$, the next three steps will go straightly. \\

{\bf STEP 1:~ Building a couple of G-FBSDEs}\\
	
By Proposition \ref{sec}, the optimal strategy $\pi^*$ is required to meet
	
	\begin{equation}\label{ssc}
		\begin{aligned}
			\dfrac{1}{\tilde{X}_T^{x_0,\pi^*}}x_0 D_T \exp & \left \{ \int_0^T \pi_s (\mu - r)\dif s + \int_0^T \pi_s (\sigma - \sigma_r)\dif B_s - \int_0^T \Big [ \dfrac{1}{2}\pi_s^2 (\sigma - \sigma_r )^2 \right. \\
			 & ~ + ~\pi_s\sigma_r (\sigma - \sigma_r) - c\pi_s \Big ] \dif \langle B \rangle_s \Bigg \} {~\text{is a symmetric G-martingale.}}
		\end{aligned}
	\end{equation}	
Set
\begin{equation}\label{Z_tilde}
	\tilde{Z}_t^*=\dfrac{\cE_t\left[\dfrac{1}{\tilde{X}_T^{x_0,\pi^*}} D_T \right]}{\cE\left[\dfrac{1}{\tilde{X}_T^{x_0,\pi^*}} D_T \right]}, \ \ t\in [0,T].
\end{equation}
Then
\begin{equation*}
	\tilde{Z}_T^* = \dfrac{\dfrac{1}{\tilde{X}_T^{x_0,\pi^*}} D_T}{\cE\left[\dfrac{1}{\tilde{X}_T^{x_0,\pi^*}} D_T \right]},\ \ \ \tilde{Z}_{\tau}^* = \dfrac{\cE_{\tau}\left[\dfrac{1}{\tilde{X}_T^{x_0,\pi^*}} D_T \right]}{\cE \left[ \dfrac{1}{\tilde{X}_T^{x_0,\pi^*}} D_T \right]},
\end{equation*}
for any stopping time $\tau \leq T$,~ q.s.. Let
$\pi_{t}^{\tau} = I_{\{t\leq \tau\}}$,~then $\pi^{\tau} \in \Pi$.

Substituting $\pi^{\tau}$ into (\ref{ssc}), one obtains
\begin{align*}
	&\cE \left[ \dfrac{1}{\tilde{X}_T^{x_0,\pi^*}} D_T \right] \cE \left[\tilde{Z}_T^* x_0 \exp \left \{\int_0^T \pi^{\tau}_s(\mu - r) \dif s + \int_0^T \pi^{\tau}_s t(\sigma - \sigma_r)\dif B_s - \right. \right.\\
	& \hspace{10em} \left.\left.\int_0^T \left[ \dfrac{1}{2}{\pi^{\tau}_s}^2 \left(\sigma - \sigma_r \right)^2 + \pi^{\tau}_s \sigma_r \left(\sigma - \sigma_r \right) - c \pi_s \right] \dif \langle B \rangle_s \right\} \right] \\
	& = \cE\left[ \dfrac{1}{\tilde{X}_T^{x_0,\pi^*}} D_T \right] \cE \left[ \cE_{\tau} \left[ \tilde{Z}_T^* \right] x_0  \exp \left \{ \!\!\int_0^{\tau} \!\!(\mu - r ) \dif s + \!\! \int_0^{\tau} \!\!(\sigma - \sigma_r) \dif B_s - \!\! \int_0^{\tau} \!\! \dfrac{1}{2} (\sigma^2 \!- \!\sigma_r^2 \!- \!2c) \dif\langle B \rangle_s \right\} \right] \!\! \\
	&= \cE \left[ \dfrac{1}{\tilde{X}_T^{x_0,\pi^*}} D_T \right] \cE \left[ \tilde{Z}_{\tau}^* x_0 \exp \left\{ \left(\mu - r \right)\tau +  (\sigma - \sigma_r) B_{\tau} - \dfrac{1}{2}(\sigma^2 - \sigma_r^2 - 2c) \langle B \rangle_{\tau} \right \} \right].
\end{align*}
It follows from ($\ref{ssc}$) that
\begin{equation*}\label{ssecmv1}
	\tilde{Z}_{t}^* x_0 \exp \left \{ (\mu - r)t + (\sigma - \sigma_r) B_{t} - \dfrac{1}{2}(\sigma^2 - \sigma_r^2 - 2c) \langle B \rangle_{t} \right \}
\end{equation*}
is a symmetric G-martingale.
Since $(\tilde{Z}_t^*)$ is a G-martingale, then by Theorem \ref{GT}, one has
\begin{equation}\label{zr_tilde}
	\dif \tilde{Z}^*_t = - \dfrac{c + \sigma_r(\sigma_r - \sigma)}{\sigma - \sigma_r} \tilde{Z}^*_t \dif B_t - \dfrac{\mu - r}{\sigma - \sigma_r} \tilde{Z}^*_t \dif \tilde{B}_t, \ \ t\in[0,T].
\end{equation}
Combining this with (\ref{X_tilde}), we conclude that the triplet $(\tilde{X}^{x_0,\pi^*},\pi^*,\tilde{Z}^*)$ is a solution of the following G-FBSDEs	
\begin{equation*}\label{sG-fbsde}
	\begin{cases}
		\dif\tilde{X}_t^{x_0,\pi} = \tilde{X}_t^{x_0,\pi} (r\dif t + \sigma_r \dif B_t) + \pi_t \tilde{X}_t^{x_0,\pi} \big ((\mu-r)\dif t + (\sigma - \sigma_r) \dif B_t + c \dif \langle B \rangle_t \big ), \ \ t \in [0,T],\\
		\tilde{X}_0^{x_0,\pi} = x_0,\\
		\dif\tilde{Z}_t = - \dfrac{c + \sigma_r(\sigma_r - \sigma)}{\sigma - \sigma_r}\tilde{Z}_t \dif B_t - \dfrac{\mu - r}{\sigma - \sigma_r} \tilde{Z}_t \dif\tilde{B}_t, \ \ t \in [0,T],\\
		\tilde{Z}_T = \dfrac{\dfrac{1}{\tilde{X}_T^{x_0,\pi}} D_T}{\cE \left [ \dfrac{1}{\tilde{X}_T^{x_0,\pi}} D_T \right ]}.
	\end{cases}
\end{equation*}	
\\
{\bf STEP 2: Solving the G-FBSDEs}\\

Firstly, for any $t \in [0,T]$, employing G-It\^o's formula,  one gets
\begin{equation*}\label{1}
	\begin{aligned}
		\dif \Big (\dfrac{1}{\tilde{X}_t^{x_0, \pi^*}} D_t  & \Big ) = ~\dfrac{1}{\tilde{X}_t^{x_0, \pi^*}} \dif D_t + D_t \dif \dfrac{1}{\tilde{X}_t^{x_0, \pi^*}} - \dfrac{D_t \sigma_r \big(\sigma_r + \pi_t^*(\sigma - \sigma_r) \big)}{\tilde{X}_t^{x_0, \pi^*}} \dif \langle B \rangle_t \\
		=& \dfrac{-\pi_t^* D_t (\mu - r ) \dif t - \pi_t^* D_t (\sigma - \sigma_r) \dif B_t + \pi_t^* D_t \big [ \big( \sigma_r + \pi_t^*(\sigma - \sigma_r) \big)(\sigma - \sigma_r) - c \big ] \dif \langle B \rangle_t}{\tilde{X}_t^{x_0,\pi^*}}.
	\end{aligned}
\end{equation*}
Define a process $(N_t)$ as follows
\begin{equation}\label{N1}
	\begin{cases}
		~\dif N_t = a_t N_t \dif t + b_t N_t \dif \langle B \rangle_t,\ \ t\in[0,T],\\
		~N_T = 1.
	\end{cases}
\end{equation}
By (\ref{N1}) and  G-It\^o's formula, one obtains
\begin{equation}\label{Gsolve}
	\begin{aligned}
		\dif \Big ( N_t \dfrac{1}{\tilde{X}_t^{x_0,\pi^*}} D_t \Big ) =& ~ N_t \dif \dfrac{1}{\tilde{X}_t^{x_0,\pi^*}}D_t + \dfrac{1}{\tilde{X}_t^{x_0,\pi^*}} D_t \dif N_t\\
		=& ~ \dfrac{D_t N_t \big[a_t - \pi_t^*(\mu - r)\big]}{\tilde{X}_t^{x_0,\pi^*}}\dif t - \dfrac{D_t N_t \pi_t^*(\sigma-\sigma_r)}{\tilde{X}_t^{x_0,\pi^*}}\dif B_t\\
		& ~ + \dfrac{D_t N_t \Big [b_t + \pi_t^* \big[ \big(\sigma_r + \pi_t^*(\sigma-\sigma_r) \big) (\sigma - \sigma_r) - c \big] \Big]}{\tilde{X}_t^{x_0,\pi^*}} \dif \langle B \rangle_t ,\ \ t \in [0,T].
	\end{aligned}
\end{equation}
To guarantee that $\Big(N_t D_t \dfrac{1}{\tilde{X}_t^{x_0,\pi^*}}\Big)$ is a symmetric G-martingale, it is sufficient that
\begin{equation}\label{G_a,b}
	\begin{cases} 
		~a_t = \pi^*_t(\mu - r),t \in [0,T],\\
		~b_t = -\pi_t^*\big[\big(\sigma_r + \pi_t^*(\sigma - \sigma_r)\big) (\sigma - \sigma_r) - c \big],\ \ t \in [0,T].
	\end{cases}
\end{equation}
Then, by (\ref{N1})$\sim$(\ref{G_a,b}), one has
\begin{equation*}
	N_t D_t\dfrac{1}{\bar{X}_t^{x_0,\pi^*}} = D_T \dfrac{1}{\bar{X}_T^{x_0,\pi^*}} + \int_t^T \dfrac{D_s N_s \pi_s^*(\sigma - \sigma_r)}{\tilde{X}_s^{x_0,\pi^*}} \dif B_s,\ \ t \in [0,T].
\end{equation*}
Observe that $\displaystyle \Big(\int_0^t \dfrac{D_s N_s \pi_s^*\ (\sigma - \sigma_r)}{\tilde{X}_s^{x_0,\pi^*}}\dif B_s \Big)$ is a symmetric G-martingale. The following holds
\begin{equation*}
	N_t D_t\dfrac{1}{\bar{X}_t^{x_0,\pi^*}} = \cE_t \left[ D_T \dfrac{1}{\bar{X}_T^{x_0,\pi^*}} \right], \ \ t \in [0,T].
\end{equation*}
This together with (\ref{Z_tilde}), one can yield
\begin{equation*}
	\dif Z_t^*=-Z_t^*\pi_t^* (\sigma-\sigma_r) \dif B_t, \ \ t \in [0,T].
\end{equation*}
Finally, comparing to (\ref{zr_tilde}), one concludes that
\begin{equation}\label{sop}
	\pi^*_t = \Big(\dfrac{(\mu - r)\rho_t + c}{\sigma - \sigma_r} - \sigma_r\Big){\Big /}\left(\sigma - \sigma_r \right),\ \ t \in [0,T].
\end{equation}
\\
{\bf STEP3: Verifying the sufficient condition}\\

Substituting (\ref{sop}) into condition ($\ref{scondition}$), an equivalent statement of condition ($\ref{scondition}$) can be established. That is, for any $\pi \in \Pi$, the process  $\Big(\dfrac{D_t F_t}{\tilde{X}_t^{x_0,\pi^*}} \Big)$ is a symmetric G-martingale, where $(F_t)$ is defined by
$$ F_t = \exp \left\{\!\!\int_0^t \! \pi_s (\mu-s) \dif s + \!\!\!\int_0^t \! \pi_s (\sigma - \sigma_r) \dif B_s - \!\!\!\int_0^t \! \big[ \dfrac{1}{2} \pi_s^2 (\sigma - \sigma_r) + \pi_s \sigma_r (\sigma - \sigma_r) - c\pi_s \big] \dif \langle B \rangle_s \right\},$$
for any $t \in [0, T]$. 

The next step is to prove the equivalency relation.

For any $\pi \in \Pi$, by G-It\^o's formula, one can derive
\begin{equation*}
	\begin{aligned}
		\dif \dfrac{D_t F_t}{\tilde{X}_t^{x_0,\pi^*}} =& ~F_t \dif \dfrac{D_t}{\tilde{X}_t^{x_0,\pi^*}} + \dfrac{D_t}{\tilde{X}_t^{x_0,\pi^*}} \dif F_t - \dfrac{D_t F_t \pi_t}{\tilde{X}_t^{x_0,\pi^*}} \Big[(\mu - s)\dif t + \big(c - \sigma_r (\sigma - \sigma_r) \big) \dif \langle B \rangle_t \Big] \\
		=&~ \dfrac{D_t F_t \pi_t}{\tilde{X}_t^{x_0,\pi^*}} \Big[ (\mu - s)\dif t + \big(c - \sigma_r (\sigma-\sigma_r)\big) \dif \langle B\rangle_t - (\mu-s) \dif t - \big(c-\sigma_r (\sigma - \sigma_r )\big) \dif \langle B \rangle_t \Big] \\
		&~+~ \dfrac{D_t F_t \pi^*_t}{\tilde{X}_t^{x_0,\pi^*}} \Big[(\mu - r)\dif t + c\dif \langle B \rangle_t - (\mu - r)\dif t - c\dif \langle B \rangle_t \Big] \\
		&~+~ \dfrac{D_t F_t \pi_t}{\tilde{X}_t^{x_0,\pi^*}}(\sigma - \sigma_r) \dif B_t - \dfrac{D_t F_t}{\tilde{X}_t^{x_0,\pi^*}(\sigma - \sigma_r)} \Big[ (\mu - r) \dif \tilde{B}_t + \big(c - \sigma_r(\sigma - \sigma_r) \big) \dif B_t \Big] \\
		=&~ \dfrac{D_t F_t}{\tilde{X}_t^{x_0,\pi^*}} \Big[ \pi_t (\sigma - \sigma_r) - \dfrac{c}{\sigma - \sigma_r} + \sigma_r \Big] \dif B_t - \dfrac{D_t F_t}{\tilde{X}_t^{x_0,\pi^*} (\sigma - \sigma_r)}(\mu - r) \dif \tilde{B}_t,\ \ t\in[0,T].
	\end{aligned}
\end{equation*}
It implies that $\pi^*_t = \Big( \dfrac{(\mu - r)\rho_t + c}{\sigma - \sigma_r} - \sigma_r \Big) \Big/ (\sigma - \sigma_r),\ \ t\in[0,T]$ meets condition ($\ref{scondition}$).	

The next theorem summarizes above analysis. 

\bt\label{st_2}
	For utility  $U(x)=\ln x$ and the stochastic interest rate model (\ref{X_tilde})~-~(\ref{saim}), an optimal strategy is $\pi^*_t = \Big( \dfrac{(\mu - r)\rho_t + c}{\sigma - \sigma_r} - \sigma_r \Big) \Big/ (\sigma - \sigma_r),\ \ t\in[0,T]$. 
\et
 Some practical implications reveal through scenario analysis of parameters. Clearly, if the bond has no volatility uncertainty, in other words, $ \sigma_r = 0 $, then the results degenerate to deterministic case. And if $ \sigma < \sigma_r $, the relationship between optimal strategy and parameter $\rho_t,~ c $ and $ \sigma - \sigma_r $ is same as the deterministic case. For fixed $ \mu,~ r,~ \rho_t ~and~ c $, there are two cases to be discussed further.
 
 (1) Either a pair of conditions $ \sigma > \sigma_r $ and $ \dfrac{(\mu - r)\rho_t + c}{\sigma - \sigma_r} < \sigma_r $, or alternative pair $ \sigma < \sigma_r $ and $ \dfrac{(\mu - r)\rho_t + c}{\sigma - \sigma_r} > \sigma_r $ hold, optimal strategy (\ref{sop}) is less than zero, i.e., in short position. 
 
 (2) Either a pair of conditions $ \sigma > \sigma_r $ and $ \dfrac{(\mu - r)\rho_t + c}{\sigma - \sigma_r} > \sigma_r $, or alternative pair $\sigma < \sigma_r$ and $\dfrac{(\mu-r)\rho_t+c}{\sigma-\sigma_r} < \sigma_r$ hold, optimal strategy (\ref{sop}) is larger than zero, i.e., in long position. 
 
 All of these explanations imply that stock position is determined by the relationship between volatility coefficients of the bond and the stock, which is also consistent with market intuition. 

\section{Conclusion}

This paper shows the feasibility of G-martingale approach to G-utility maximization problems. For a continuous-time financial market model under uncertainty characterized by G-Brownian motion, an optimal investment strategy can be derived by the approach. Surely it depends on the special case of log-utility function, in which representative investor is myopic and focuses on logarithmic wealth level. Similar to \cite{Wxz}, a sufficient condition to solve the model is presented. However, it is difficult to give an explicit solution because the quadratic variation process $(\langle B \rangle_t)$ of G-Brownian motion is a stochastic process, but not deterministic function $t$ in the Brownian motion situation. Particularly, $\langle B \rangle_1$ is a maximal distribution on $[\underline{\sigma}^2, \overline{\sigma}^2]$. Thus, there exists no constant $c$ such that $(\langle B \rangle_t-ct)$ is a symmetric G-martingale.  A feasible G-martingale approach is carried out in an extended nonlinear expectation space.As an application of these arguments, we reconsider the stochastic interest model under uncertainty, and obtain an alternative powerful formulation and plentiful implications.  

There are two directions deserving to pay attention to the further application of nonlinear expectation theory. One is the G-martingale approach for option pricing in mathematical finance. The other one is the optimal growth under uncertainty in macroeconomics. It has appeared that the broadness of nonlinear expectation theory may be more suitable for the characterization of macroeconomic issues, with confusing uncertainty from short cycles to the long trend in decades.\\


\begin{thebibliography}{00}\addtolength{\itemsep}{0.5ex}
\bibitem{alp95}
Avellaneda, M., Levy, A., Paras, A., 1995. Pricing and Hedging Derivative Securities in Markets with Uncertain Volatilities. Applied Mathematical Finance 2, 73-88.	

\bibitem{be19}
Beissner, P., 2019. Coherent-Price System and Uncertainty-Neutral Valuation. Risk 98 (7), 1-18.

\bibitem{Bmb}
Biagini, F., Mancin, J., Brandis, T., 2019. Robust Mean–Variance Hedging via G-Expectation. Stochastic Processes and Their Applications 129 (4), 1287-1325.

\bibitem{Ce}
Chen Z., Epstein, L., 2002. Ambiguity, Risk and Asset Returns in Continuous Time. Econometrica 70 (4), 1403-1443.

\bibitem{Cho}
Choquet, G., 1953. Theory of Capacity. Annales de l’Institute Fourier(Grenoble) 5, 31-295.

\bibitem{Do}
Dothan, U., 1978. On the Term Structure of Interest Rate. Journal of Financial Economics 6, 59-69.

\bibitem{Ej13}
Epstein, L., Ji, S., 2013. Ambiguous Volatility and Asset Pricing in Continuous Time. The Review of Financial Studies 26, 1740-1786.

\bibitem{Ej14}
Epstein, L., Ji, S., 2014. Ambiguous Volatility, Possibility and Utility in Continuous Time. Journal of Mathematical Economics 50, 269-282.

\bibitem{GS}
Gilboa, I., Schmeidlerand, D., 1989. Maxmin Expected Utility with Non-Unique Prior. Journal of Mathematical Economics 18, 141-153.

\bibitem{Hj14}
Hu, M., Ji, S., Peng, S., Song, Y., 2014. Comparison Theorem, Feynman–Kac Formula and Girsanov Transformation for BSDEs Driven by G-Brownian Motion.  Stochastic Processes and Their Applications 124, 1170–1195.

\bibitem{Hj18}
Hu, M., Ji, S., 2018. A Note on Pricing of Contingent Claims under G-Expectation. ArXiv:1303.4274V1.

\bibitem{ho2101}
 H$\ddot {\rm o}$lzermann, J., 2021. The Hull-White Model under Volatility Uncertainty. ArXiv:1808.03463v3.

\bibitem{ho2102}
 H$\ddot {\rm o}$lzermann, J., 2021. Term Structure Modeling under Volatility Uncertainty. ArXiv:1904.02930v4.

\bibitem{ho2103}
 H$\ddot {\rm o}$lzermann, J., 2021. Pricing Interest Rate Derivatives under Volatility
Uncertainty. ArXiv:2003.04606v3.

\bibitem{Keynes}
Keynes, J., 1921. A Treatise of Probability. Science 58, 51-52.

\bibitem{Kni}
Knight, F., 1921. Risk, Uncertainty and Profit. Social Science Electronic Publishing (4), 682-690.

\bibitem{Lr}
Lin, Q., Riedel, F., 2021. Optimal Consumption and Portfolio Choice with Ambiguous Interest Rates and Volatility. Economic Theory 71, 1189-1202.

\bibitem{Lyo}
Lyons, T., 1995. Uncertain Volatility and the Risk-Free Synthesis of Derivatives. Applied Mathematical Finance 2, 7–133.

\bibitem{PLyo}
Peng, S., 2004. Filtration Consistent Nonlinear Expectations and Evaluations of Contingent Claims. Acta Mathematicae Applicatae Sinica 20, 1-24.

\bibitem{Pe97}
Peng, S., 1997. Backward SDE and Related g-Expectation, in Backward Stochastic Differential Equations, Pitman Research Notes in Math. Series 364,  141-159.

\bibitem{Pe04}
Peng, S., 2004. Filtration Consistent Nonlinear Expectations and Evaluations of Contingent Claims.  Acta Mathematicae Applicatae Sinica 20, 1-24.

\bibitem{Pe05}
Peng, S., 2005. Nonlinear Expectation and Nonlinear Markov Chains. Chinese Annals of Mathematics B 26 (2), 159-184.

\bibitem{Pe06}
Peng, S., 2006. G-Expectation, G-Brownian Motion and Related Stochastic Calculus of It\^o's Type. Stochastic Analysis and Applications Abel Symposia 2, 541-567.

\bibitem{Pe17}
Peng, S., 2017. Theory, Methods and Meaning of Nonlinear Expectation Theory (in Chinese). Sci Sin Math 47, 1223–1254.

\bibitem{Pe19}
Peng, S., 2019. Nonlinear Expectations and Stochastic Calculus under Uncertainty. Springer.

\bibitem{Pe14}
Peng, S., Song, Y., Zhang, J., 2014. A Complete Representation Theorem for G-martingales. Stochastics 86 (4), 609-631.

\bibitem{pyy21}
Peng, S., Yang, S., Yao, J., 2021. Improving Value-at-Risk Prediction under Model Uncertainty. Journal of Financial Econometrics, 1-32.

\bibitem{py21}
Peng, S., Yang, S., 2021. Autoregressive Models of Time Series under Volatility Uncertainty and Application to VaR Model. ArXiv:2011.09226v1.


\bibitem{Sz}
Sun, Z., 2018. Upper Bounds for Ruin Probabilities under Model Uncertainty. Communications in Statistics - Theory and Methods 18, 4511-4527.

\bibitem{X}
Xiao, X., 2014. Stochastic Dominance under the Nonlinear Expected Utilities. Mathematical Problems in Engineering, 1-6.

\bibitem{V}
Vorbrink, J., 2014. Financial Markets with Volatility Uncertainty. Journal of Mathematical Economics 53, 64-78.

\bibitem{Wan}
Wang, S., 2000. Class of Distortion Operators for Pricing Financial and Insurance Risks. Journal of Risk and Insurance 7 (1), 15-36.

\bibitem{Wxz}
Wang, Z., Xia, J., Zhang, L., 2007. Optimal Investment for an Insurer: The Martingale Approach. Insurance: Mathematics and Economics 40, 322-334.

\end{thebibliography}
\end{document}